\documentclass[12pt]{article}
\usepackage{amsmath}
\usepackage{amsthm}
\usepackage{amssymb}
\usepackage{graphicx}
\theoremstyle{plain}
\newtheorem{theorem}{Theorem}[section]

\newtheorem{lemma}[theorem]{Lemma}

\newtheorem{proposition}[theorem]{Proposition}

\newtheorem{corollary}[theorem]{Corollary}

\newtheorem{remark}[theorem]{Remark}

\newtheorem{example}[theorem]{Example}

\theoremstyle{definition}
\newtheorem{definition}[theorem]{Definition}

\begin{document}
\newcommand{\field}[1]{\ensuremath{\mathbb{#1}}}
\newcommand{\N}{\field{N}}
\newcommand{\Z}{\field{Z}}

\title{The subword complexity of a class of infinite binary words}

\author{Irina\ Gheorghiciuc \thanks{Mathematical Sciences Department, 
University of Delaware, Newark, DE 19716--2553, {\tt <gheorghi@math.udel.edu>}. 
This paper represents a modified portion of the author's PhD dissertation at the University of Pennsylvania, thesis adviser Herbert S. Wilf.} } 

\maketitle

\begin{abstract}

Let $A_q$ be a $q$-letter alphabet and $w$ be a right infinite word on this alphabet. A subword of $w$ is a block of consecutive letters of $w$. The subword 
complexity function of $w$ assigns to each positive integer $n$ the number $f_w(n)$ 
of distinct subwords of length $n$ of $w$.

The gap function of an infinite word over the binary alphabet $\{ 0,1 \}$ gives the distances between consecutive $1$'s in this word. In this paper we study infinite binary words whose gap function is injective or "almost injective". A method for computing the subword complexity of such words is given. A necessary and sufficient condition for a function to be the subword complexity function of a binary word whose gap function is strictly increasing is obtained.

\end{abstract}

\section{Introduction}
\label{sec:int}

The subword complexity (sometimes called symbolic complexity) of finite and infinite words became an important subject in Combinatorics on Words recently. Applications include Dynamical Systems, Ergodic Theory and Theoretical Computer Science (\cite{Ferenczi1999}).  

For a given infinite word it is not easy to compute the subword complexity function. Classes of infinite words whose subword complexity function
has been computed include paperfolding sequences (see \cite{Allouche1992}), Rudin-Shapiro sequences, Thue-Morse sequences and generalized
Thue-Morse sequences (see \cite{TrompShallit}) and sequences defined by billiards
in hypercubes, which generalize Sturmian words. A survey of results of this kind can be found in \cite{Allouche1994} and \cite{AllucheShallit}, with \cite{AllucheShallit} being the most recent. 

Another general problem of much interest is to determine which function can be the subword complexity function of an infinite word. A list of known necessary conditions as well as a list of sufficient conditions is given in \cite{Ferenczi1999}. 

If $u$ and $v$ are two finite words over the same alphabet, then $uv$ will denote the concatenation of $u$ and $v$. In particular, for a positive integer $n$, $u^n=uu...u$ ($n$ times). $u^0= \epsilon$, where $\epsilon$ is the empty word.
 
Let $w$ be an infinite word over the binary alphabet $A=\{ 0,1 \}$. 
Since the subword complexity of a binary word does not change 
when we interchange $0$ and $1$, we can assume without loss 
of generality that $w$ contains an infinite number of $1$'s.

\begin{definition}
The function $G:\N \rightarrow \N $ is called the $1$-distribution function of 
$w$ if $G(i)$ is the position of the $i$th $1$ in $w$. By convention $G(0)=0$. 
The function $g(i)=G(i)-G(i-1)$ defined for $i \geq 1$ is called the gap
function. If $g$ is strictly increasing, then $w$ is said to be gap increasing.
\end{definition}

The $1$-distribution function $G(i)$ is called sometimes the occurrence function of a letter (see \cite{AllucheShallit}). The main result of Section~\ref{sec:gapinc} is the following theorem.

\begin{theorem} [The Subword Complexity of Gap Increasing Words] Let $w$ be a
gap increasing word with $1$-distribution function $G$ and gap function $g$.
If $n \leq g(1)$, then the subword  complexity function of $w$ is $f_w(n)=n+1$,
otherwise  
\begin{equation*}
f_w(n)=G(L_n)+G(L_n+1)-G(M_n+1)+n(M_n-L_n)+n+1,
\end{equation*} 
where  $L_n$ is the least non-negative integer which satisfies
\begin{equation*}
g(L_n+1)+g(L_n+2)\geq n+1 
\end{equation*} 
and $M_n$ is the maximum non-negative integer which satisfies 
\begin{equation*}
g(M_n+1)\leq n-1.
\end{equation*}
\label{gapinc} 
\end{theorem}

Given a function $g$ that satisfies the hypothesis of Theorem~\ref{gapinc}, we always can compute $M_n$ and $L_n$ in $O(log(n))$ time.

We use Theorem~\ref{gapinc} to find the asymptotic behavior of the
subword complexity function of infinite words over $A=\{ 0,1 \}$ with
polynomial and exponential $1$-distribution functions. Let $k>1$ be an
integer. If the $1$-distribution function of $w$ is $G(n)=n^{k}$, then
$f_w(n)=\Theta (n^{\frac {k}{k-1}})$. If the $1$-distribution function of $w$
is $G(n)=k^{n-1}$, then $f_w(n)=\Theta (n)$.

In Section~\ref{sec:rec} we give a recurrence formula for the subword
complexity function  of infinite binary words whose gap function is injective
and later, in Section~\ref{sec:block}, generalize it to binary words whose gap function is blockwise injective. A function  $h: \N \rightarrow \N$ is called blockwise injective if $h(i)=h(j)$ for $i<j$ implies $h(i)=h(i+1)=...=h(j)$.
We prove that if $w$ has the gap function $g(n)=cn+d$, where $c,d\in\N$ and $c\geq 2$, then $f_w(n+2c)=f_w(n)+n+2c$ for $n \geq c+d$.

In Section~\ref{sec:cond} we derive necessary and sufficient conditions for a
function  $f:\N \rightarrow \N$ to be the subword complexity of an infinite
gap  increasing word. First we prove the following necessary condition: $f(1)=2$, $f(2)=3$ and $| { \Delta }^2 f(n) | \leq 1$ ($\Delta$ is the usual difference operator). 

To formulate the necessary and sufficient condition, we introduce additional notation. Let $\{ x_i  \} _{i \geq 1}$ and $\{ y_i  \} _{i \geq 1}$ be two strictly increasing sequences that have no elements in common, then $\{ x_i  \} _{i \geq 1} \sqcup \{ y_i  \} _{i \geq 1}$ will denote the strictly increasing sequence that consists of all the elements of both sequences. We allow a sequence to have no elements (be empty). If $\{ y_i  \} _{i \geq 1}$ is empty, then $\{ x_i  \} _{i \geq 1} \sqcup \{ y_i  \} _{i \geq 1} = \{ x_i  \} _{i \geq 1}$. Also, if $\{ x_i  \} _{i \geq 1}$ is a strictly increasing sequence of positive integers and $p<x_1$ be a non-negative integer, then $\sigma_p ( \{ x_i  \} _{i \geq 1} )$ will denote the strictly increasing sequence $\{ p+x_1+1,x_1+x_2+1,x_2+x_3+1,... \} $.

\begin{theorem} [Necessary and sufficient condition for a function to be the
subword complexity function of an infinite gap increasing word] Let $f: \Z _+
\rightarrow \Z _+$ be such that $f(1)=2$, $f(2)=3$ and $| { \Delta }^2 f(n) | \leq 1$ (this is a necessary condition). Set $\{ a_i  \} _{i \geq 1}$ to the sequence obtained by arranging the elements of the set $\{ n \in \Z _+ |\Delta ^2 f(n) = 1 \}$ in strictly increasing order. 
Similarly $\{ b_i  \} _{i \geq 1}$ is obtained by arranging the elements of $\{ n \in \Z _+ | \Delta ^2 f(n) = -1 \}$ in strictly increasing order. Then $f$ is the subword complexity function of an infinite gap
increasing word if and only if $\{ a_i  \} _{i \geq 1}$ is not empty and there exist an integer $p$, 
$0 \leq p < a_1$, and a strictly increasing sequence of positive integers $\{ c_i  \} _{i \geq 1}$, 
with no elements in common with either $\{ a_i  \} _{i \geq 1}$ or $\{ b_i  \} _{i \geq 1}$, 
such that $\{ a_i  \} _{i \geq 1} \sqcup \{ c_i  \} _{i \geq 1}$ is infinite and 
\begin{equation}
\sigma_p (\{ a_i  \} _{i \geq 1} \sqcup \{ c_i  \} _{i \geq 1}) = 
\{ b_i  \} _{i \geq 1} \sqcup \{ c_i  \} _{i \geq 1}.
\end{equation}
\label{th:prenecsuf}   
\end{theorem} 

Based on Theorem~\ref{th:prenecsuf}, in Section~\ref{sec:cond} we describe a practical method to determine if a function $f$ is the subword complexity function of an infinite gap increasing word and give examples (see Examples~\ref{ex:yes} and \ref{ex:no}). We also show that it is not true in general that a function $f:N \rightarrow N$ cannot be the subword complexity function of two distinct infinite gap increasing words. 

Many sequences that have been studied thoroughly (like Sturmian, Arnoux-Rauzy and Rote and Rudin-Shapiro sequences) have affine or ultimately affine subword complexity functions (see \cite{Morse1938}, \cite{Morse1940}, \cite{ArnouxRauzy}, \cite{Rote}). An infinite word $w$ has subword complexity $\phi (n)$  ultimately if there exists a positive integer $N$, such that for all integers $n \geq N$, the subword complexity function of $w$ is $\phi (n)$.

Cassaigne proved in \cite{Cassaigne} that if $a$ and $b$ are two integers, then there exists an infinite binary word with subword complexity function $an+b$ ultimately if and only if one of the two following conditions holds: 

1) $a \geq 2$;
 
2) $0 \leq a \leq 1$ and $b \geq 1$. 

In Section~\ref{sec:aff} we prove that for two integers $a$ and $b$ there exists a gap increasing word $w$ with the subword complexity function $f_w(n)=an+b$ ultimately if and only if $a \geq 2$. This result is generalized in the next proposition, proved in Section~\ref{sec:block}. 

\begin{proposition} Let $a$ and $b$ be two integers. There exists an infinite binary word $w$ with a blockwise injective gap function and subword complexity function $f_w(n)=an+b$ ultimately if and only if one of the two following conditions holds:

1) $a \geq 2$. 

2) $a=0$ and $b \geq 1$.
\end{proposition}

In Section~\ref{sec:aff}, in order to investigate when $an+b$ is ultimately the subword complexity of a gap increasing word, we show a geometric interpretation 
of the subword complexity of a gap increasing word. A consequence is that if $w$ is gap increasing, then the subword complexity $f$ of $w$ satisfies $n+1 \leq f(n+1) \leq \left \lceil n/2 \right \rceil \left \lfloor n/2 \right \rfloor +\left \lceil n/2 \right \rceil +1$ and both bounds are exact.

The main result of Section~\ref{sec:block} is Proposition~\ref{pr:bl}, which gives a method of computing the subword complexity of infinite binary words whose gap function is blockwise injective. A consequence of Proposition~\ref{pr:bl} is that, if $w$ is an infinite binary word whose gap function is blockwise injective, then its subword complexity $f_w(n)=O(n^3)$, more precisely $1 \leq f(n) \leq \frac {n^3}{6} + \frac {5n}{6} +1$. Thus the topological entropy of $w$ is $0$.

Since Proposition~\ref{pr:bl} is rather big for an introduction section, here we present a method of computing the subword complexity function of infinite binary words whose gap function is non-decreasing (a non-decreasing function is blockwise injective). 

To each unbounded non-decreasing function $g:\N \rightarrow \N$ we assign two
functions $j:\N \rightarrow \N$ and $p:\N \rightarrow \N$ such that $j$ is the
strictly increasing function which assumes the same values and in the same
order as $g$ does, and, for each $r \in \N$, $p(r)$ is the number of times $g$
assumes value $j(r)$.  If $g:\N \rightarrow \N$ is a bounded
non-decreasing function, let $b$ be the  number of distinct values that $g$
assumes. Then  $j(r)$ is defined as above for  $1 \leq r \leq b$, and $p(r)$
is defined for $1 \leq r \leq b-1$. 

\begin{lemma} Let $w$ be an infinite binary word whose gap function $g$
is non-decreasing and $j(r)$, $p(r)$ (and $b$ if $g$ is bounded) be 
defined as above.

If $g$ is unbounded, then, for every natural $n$, $f_w(n+1)-f_w(n)$ equals 
one plus the number of integer $r$ solutions of   
\begin{equation}    
    \begin{cases} 
    j(r)\leq n \\
    j(r-1)+p(r)j(r)\geq n+1 .
    \end{cases}
    \label{eq:nondec1'}
\end{equation}
If $g$ is bounded, for every natural number $n$ consider the inequality system 
\begin{equation}
    \begin{cases} 
    r<b \\
    j(r)\leq n \\
    j(r-1)+p(r)j(r)\geq n+1 .
    \end{cases}
    \label{eq:nondec2'}
\end{equation}
If $1 \leq n \leq j(b)-1$, then $f_w(n+1)-f_w(n)$ equals one plus the number of 
integer $r$ solutions of eq.~(\ref{eq:nondec2'}). If $n \geq j(b)$, then
$f_w(n+1)-f_w(n)$ is  just the number of integer solutions of eq.~(\ref{eq:nondec2'}).
\label{lm:non'}
\end{lemma}

\section{Gap increasing words}
\label{sec:gapinc}

An infinite word is said to be gap increasing if its gap function is strictly 
increasing. The main result of this section is the following theorem. 

\begin{theorem} Let $w$ be a gap increasing word with $1$-distribution
function $G$ and gap function $g$. For $n \leq g(1)$ the subword complexity function of $w$ is $f_w(n)=n+1$. For $n>g(1)$ the subword complexity function of $w$ is  
\begin{equation}
f_w(n)=G(L_n)+G(L_n+1)-G(M_n+1)+n(M_n-L_n)+n+1,
\end{equation} 
where  $L_n$ is the least non-negative integer which satisfies
\begin{equation}
g(L_n+1)+g(L_n+2)\geq n+1 
\label{eq:lnin}
\end{equation} 
and $M_n$ is the maximum non-negative integer which satisfies 
\begin{equation}
g(M_n+1)\leq n-1.
\label{eq:mnin}
\end{equation} 
\label{th:gapinc} 
\end{theorem}

The proof of Theorem~\ref{th:gapinc} follows from Propositions~\ref{pr:n+1} and 
\ref{pr:n+1} proved later in this section. 

It should be mentioned that, given a function $g$ that satisfies the hypothesis
of Theorem~\ref{th:gapinc} and $n > g(1)$, we always can compute $M_n$ and $L_n$ 
in $log(n)$ time. To see this notice that, since $g$ is strictly increasing, $g(i+1) 
\geq i$ for all $i \geq 0$ and therefore $g(i+1)+g(i+2) \geq 2i+1$ for all $i \geq 0$. 
This implies that $L_n \leq n$ and $M_n < n$. We can use the dichotomy algorithm 
with initial value $n$ for both $L_n$ and $M_n$ to find $L_n$ as the minimum 
solution of eq.~(\ref{eq:lnin}) and $M_n$ as the maximum solution of 
eq.~(\ref{eq:lnin}) (here we use again the fact $g$ is strictly increasing). 
The running time for both is $log(n)$.  

\begin{proposition} Let $w$ be an infinite gap increasing word. The number of 
subwords of length $n$ of $w$ that contain at most one $1$ is $n+1$. 
\label{pr:n+1}
\end{proposition}

\begin{proof} Let $g$ be the gap function of $w$. Since $g$ is strictly increasing, 
we can pick $L$ such that $g(L)\geq n+1$. Then $w'=10^{g(L)-1}10^{g(L+1)-1}1$ is a subword $w$ which contains
all possible words of length $n$ over $\{ 0,1 \}$ that have at most one $1$. 
There are $n+1$ such words, thus the number of subwords of $w$ of length $n$ 
that contain at most one $1$ is $n+1$.
\end{proof}

\begin{proposition} Let $w$ be an infinite gap increasing word with gap function $g$ 
and $1$-distribution function $G$. The number of subwords of $w$ of length $n$ that 
contain at least two $1$'s is $G(L_n)+G(L_n+1)-G(M_n+1)+n(M_n-L_n)$ if $n > g(1)$, 
where $L_n$ and $M_n$ are defined in Theorem~\ref{th:gapinc}, and $0$ if $n \leq g(1)$.
\label{pr:M_n}
\end{proposition}

\begin{proof} Fix a positive integer $n$. Let $v_i$ be the subword of length $n$ of $w$ which occurs at place 
$i$ in $w$ and $S$ be the number of subwords of $w$ of length $n$ that contain 
at least two $1$'s. 

Since $g$ is strictly increasing, all subwords of $w$ which contain at least two 
$1$'s, occur just once in $w$. Thus $S$ equals the number of $i$'s for which $v_i$ 
contains at least two $1$'s.

Let $N$ be the least integer such that $v_N$ contains at most one $1$.
Then the number of subwords of $w$ of length $n$ that contain at least two $1$'s
$S=N-1+K$, where $K$ is the number of $i$'s such that $v_i$ 
contains at least two $1$'s and $i>N$.

First we find N. There are two possibilities: either $N=1$, which happens
when $n \leq g(1)$ and then $S=0$, or the $(N-1)$th character of $w$ 
is $1$ (otherwise $v_{N-1}$ would contain at most one $1$, which would contradict 
the choice of $N$). In the last case $N-1=G(L_n)$ for some positive integer $L_n$.
  
\begin{equation}  
w=\dots
01\overbrace{ \underbrace {00\dots 00}_{g(L_n+1)-1} 100\dots 0}^{v_N}0\dots 010\dots  
\end{equation} 

Because $N=G(L_n)+1$ is the least integer for which $v_N$ contains at most one $1$,
we have that $L_n$ is the least integer which satisfies the inequality 
$g(L_n+1)+g(L_n+2)-1 \geq n$. 

Next we find the number $K$ of $i$'s such that $i>G(L_n)$ and $v_i$ contains 
at least two $1$'s. All such $v_i$'s contain exactly two $1$'s (otherwise there
would exist an $l>L_n$ such that $g(l+1)+g(l+2)+1 \leq n$ which is impossible 
because $g(l+1)+g(l+2)+1 > g(L_n+1)+g(L_n+2)+1 \geq n$).

Let $P$ be the maximum integer for which $v_P$ contains exactly two $1$'s.
Clearly the first character of $v_P$ is $1$, otherwise $v_{P+1}$ 
would contain two $1$'s too, which would contradict the choice of $P$. So
$P=G(M_n)$ where $M_n$ is the maximum integer which satisfies
$g(M_n+1)+1 \leq n$.

Hence $K$ is the number of $v_i$'s, $G(L_n)<i\leq G(M_n)$,
that contain exactly two $1$'s. For each integer $l$, $L_n<l\leq M_n$, we
compute the number of $i$'s in the interval $G(l-1)<i\leq G(l)$ for which
$v_i$ contains two $1$'s. If $G(l-1)<i\leq G(l)$, then, as shown in  (\ref{pic}), 
the only two $1$'s that $v_i$ can contain are the $G(l)$th and $G(l+1)$th 
characters of $w$.

\begin{equation}
w=\dots 0\underset{G(l-1)}{1}0\dots 00\overbrace{0\dots 00 \overbrace
{\underset{G(l)}{1}00\dots 00\underset{G(l+1)}{1}}^{g(l+1)+1}00\dots
0}^{v_i} 0\dots 0\underset{G(l+2)}{1}0\dots
\label{pic}
\end{equation}

Thus $v_i$ contains two $1$'s if and only if $G(l+1)-n+1 \leq i \leq G(l)$ and
there are $G(l)-G(l+1)+n$ such $v_i$'s.

The total number of $v_i$'s which contain two $1$'s and such that $i>G(L_n)$
is $$K = \sum_{l=L_n+1}^{M_n} (G(l)-G(l+1)+n)=G(L_n+1)-G(M_n+1)+n(M_n-L_n).$$
Thus the number of subwords of $w$ of length $n$ that contain at least two $1$'s
is $$S=N-1+K=G(L_n)+G(L_n+1)-G(M_n+1)+n(M_n-L_n),$$ which is the claim of the
Proposition.
\end{proof}

Next we use Theorem~\ref{th:gapinc} to find the asymptotic behavior of the
subword complexity function of infinite words with polynomial and exponential
$1$-distribution functions.

\begin{lemma} Let $k>1$ be an integer and $w$ be an infinite word whose
$1$-distribution function is $G(n)=n^{k}$. The subword complexity function of
$w$ is $$f_w(n)=\Theta (n^{\frac {k}{k-1}}),$$ where the $\Theta$-notation depends on $k$. 
\end{lemma}

\begin{proof} The gap function of $w$ is $g(n)=n^{k}-(n-1)^{k}$ and 
eq.~(\ref{eq:lnin}) becomes $(l+2)^{k}-l^{k}\geq n+1$. If $n > k2^k$, 
then any $l$ which satisfies eq.~(\ref{eq:lnin}) is greater than 2, 
thus ${(l+2)^{k} -l^{k} < 2k} {k \choose \left \lfloor \frac {k}{2}\right \rfloor} 
l^{k-1} = c_1 l^{k-1}$  which yields $L_n > (\frac {n+1}{c_1})^{1/(k-1)}$. 
On the other hand $(l+2)^{k} -l^{k} > 2kl^{k-1}$,  so $L_n < (\frac {n+1}{2k})^{1/(k-1)}$. Hence
$L_n = \Theta (n^{\frac {1}{k-1}})$,  where the $\Theta$-notation depends on
$k$. In the same way $M_n = \Theta (n^{\frac {1}{k-1}})$.  By Theorem~\ref{th:gapinc}
$f_w(n)=\Theta (n^{\frac {k}{k-1}})$.  
\end{proof}

\begin{lemma} For an integer $k>1$ consider the infinite word with
$1$-distribution function $G(n)=k^{n-1}$. The subword complexity function of $w$  
$$f_w(n)=\Theta (n),$$ where the $\Theta$-notation does not depend on $k$. 
\end{lemma}

\begin{proof} Solving eq.~(\ref{eq:lnin}) and  (\ref{eq:mnin}) for $L_n$ and $M_n$ 
respectively we get $L_n=1+\left \lceil \log _k  \frac
{n+1}{k^2-1}\right \rceil$ and $M_n=1+\left \lfloor \log _k  \frac
{n-1}{k-1}\right \rfloor$. By Theorem~\ref{th:gapinc} $f_w(n)=$ 
$$k^{\left \lceil \log _k  \frac{n+1}{k^2-1}\right \rceil}+
k^{\left \lceil \log _k  \frac{n+1}{k^2-1}\right \rceil +1}-k^{\left \lfloor \log _k  \frac
{n-1}{k-1}\right \rfloor +1}+n(\left \lfloor \log _k  \frac
{n-1}{k-1}\right \rfloor - \left \lceil \log _k  \frac{n+1}{k^2-1}\right \rceil)+n+1.$$ 
As to the asymptotic behavior, $f_w(n)=\Theta (n)$.
\end{proof}

\section{Recurrence formulas}
\label{sec:rec}   

The formula that we obtained in Theorem~\ref{th:gapinc} can be used to deduce a recurrence for the subword complexity function of an infinite gap increasing word. However a combinatorial approach that we present next leads to a
more elegant form of recurrence  and applies to a larger class of infinite
binary words, those whose gap function is injective. In Section~\ref{sec:block} the method of this section is generalized to infinite binary words whose gap function is blockwise injective (the definition of blockwise injectivity will be given later, in Section~\ref{sec:block}; an instance of blockwise injective functions are the non-decreasing functions). 

The following definition applies to finite as well as infinite binary words. 

\begin{definition}
A subword $u$ of a binary word $w$ is called a right (left) special factor of $w$ if 
both $u0$ and $u1$ ($0u$ and $1u$ respectively) are subwords of $w$. Let
$s_w(n)$  denote the number of right special factors of $w$ of length $n$
and $s'_w(n)$  denote the set of all left special factors of $w$ of length
$n$. 
\end{definition}

A very general recurrence formula for the subword complexity of an infinite binary 
word is
\begin{equation}
f_w(n+1)=f_w(n)+s_w(n)
\label{eq:spec}
\end{equation}

In Proposition~\ref{pr:inj} we show that the number $s_w(n)$ of right special
factors of length n of an infinite binary word $w$ whose gap function $g$ is
injective equals the number of solutions of a certain system of inequalities
involving $g$ and $n$. In Proposition~\ref{pr:bl} we do
the same for infinite binary words whose gap function is blockwise injective. 

\begin{proposition} Let $w$ be an infinite binary word whose gap function $g$
is injective.  Then the number $s_w(n)$ of right special factors of $w$ of length $n$ equals {\bf one} plus the number of integer $l$ solutions of 
\begin{equation}
   \begin{cases} 
    g(l)\leq n \\
    g(l-1)+g(l)\geq n+1 .
    \end{cases}
    \label{eq:inj}
\end{equation}
\label{pr:inj}
\end{proposition}

\begin{proof} Since $g$ is injective, any subword of $w$ which contains at
least two $1$'s occurs only once in $w$, thus this subword cannot be a  right
special factor. This implies that all right special factors of $w$ contain at
most one $1$. Let $u$ be a subword of $w$ of length $n$ which contains at most
one $1$. Since $u0$ is always a subword of $w$, for $u$ to  be a right
special factor it is enough that $u1$ be a subword of $w$. Thus a subword  $u$
of $w$ is a right special factor of $w$ if and only if it contains  at
most one $1$ and $u1$ is a subword of $w$. Hence $s_w(n)$ equals the number
of subwords of $w$ of length $n+1$ which contain at most two $1$'s and whose
last letter is $1$. Obviously $0^{n}1$ is such a word and is the only such word
which contains less then two $1$'s. Thus $s_w(n)$ is the number of subwords of
$w$ of length $n+1$, which contain exactly  two $1$'s and whose last letter is
$1$, plus one.

We shall count the subwords of $w$ of length $n+1$ of form $v=0^{x}10^{g(l)-1}1$,
where $x=n-g(l)$ and $0 \leq x \leq g(l-1)-1$.
Since each such $v$ occurs just once in $w$ (because it contains two $1$'s),
we count the number of $l$'s which satisfy the inequality $0 \leq
n-g(l)\leq g(l-1)-1$, that is the number of $l$'s which satisfy eq.~(\ref{eq:inj}). 
\end{proof}

\begin{corollary} Let $w$ be a gap increasing word and $L_n$ and $M_n$ be
defined as in Theorem~\ref{th:gapinc}. Then the recurrence formula for the
subword complexity function of $w$ when $n \geq g(1)$ is 
$$f_w(n+1)=f_w(n)+M_{n+1}-L_n+1.$$ 
For the values of $n$ that are less than $g(1)$, the 
recurrence is $f_w(n+1)=f_w(n)+1$.  
\label{cr:inc}
\end{corollary}

\begin{proof} By Proposition~\ref{pr:inj} the number $s_w(n)$ of right special 
factors of $w$ of length $n$ equals one plus the number of solutions 
of eq.~(\ref{eq:inj}).

If $n < g(1)$, the number of solutions of eq.~(\ref{eq:inj})
is zero. Thus $s_w(n)=1$ and $f_w(n+1)=f_w(n)+1$. 
  
Now consider $n \geq g(1)$. We make use of the fact that g is strictly increasing. 
$g(l)\leq n$ if and only if $l\leq M_{n+1}+1$, also $g(l-1)+g(l)\geq n+1$ if and 
only if $l\geq L_n+2$. So there are $M_{n+1}-L_n$ $l$'s which satisfy eq.~(\ref{eq:inj}).
Thus, by Proposition~\ref{pr:inj}, $s_w(n)=M_{n+1}-L_n+1$ and $f_w(n+1)=f_w(n)+M_{n+1}-L_n+1$.  
\end{proof}
  
\begin{example} Consider the infinite binary word $w$ whose $1$-distribution
function is $G(n)=n^{2}$. The gap function of $w$ is $g(n)=G(n)-G(n-1)=2n-1$.
Since $g$ is strictly increasing, we can use Corollary~\ref{cr:inc} to compute 
the subword complexity of $w$. The recurrence $f_w(n+1)=f_w(n)+M_{n+1}-L_n+1$ 
starts with $n=1$ because $g(1)=1$.
 
$L_n$ is the least $l$ which satisfies $2l+1+2l+3\geq n+1$, thus $L_n=\left
\lceil (n-3)/4 \right \rceil$. $M_n$ is the maximum $m$ which satisfies
$2m+1\leq n-1$, thus $M_n=\left \lfloor n/2 \right \rfloor -1$.
By Corollary~\ref{cr:inc} $$f_w(n+1)=f_w(n)+\left \lfloor (n+1)/2 \right \rfloor -\left
\lceil (n-3)/4 \right \rceil ,$$ 
$$f_w(n+4)-f_w(n)=M_{n+1}+M_{n+2}+M_{n+3}+M_{n+4}-L_{n}+L_{n+1}+L_{n+2}+L_{n+3}+4=$$ 
$$ 
\left \lfloor (n+1)/2 \right \rfloor +\left \lfloor (n+2)/2 
\right \rfloor +\left
\lfloor (n+3)/2 \right \rfloor +\left \lfloor (n+4)/2 \right \rfloor -$$ $$\left
\lceil (n-3)/4 \right \rceil -\left \lceil (n-2)/4 \right \rceil -\left \lceil
(n-1)/4 \right \rceil -\left \lceil n/4 \right \rceil =n+1+n+3-n=n+4.$$
We get an elegant recurrence for the subword complexity function of 
$w$: $$f_w(n+4)=f_w(n)+n+4.$$
\end{example}

Next we generalize this result to all infinite words whose gap function is linear.

\begin{lemma} Let $c>0$ and $d$ be two integers such that $c+d>0$. The subword
complexity function of the infinite binary word $w$ with gap function
$g(n)=cn+d$ satisfies the recurrence $f_w(n+2c)=f_w(n)+n+2c$ for $n \geq c+d$.
\label{lm:cn+d}
\end{lemma}

\begin{proof} After solving eq.~(\ref{eq:lnin}) and (\ref{eq:mnin}), we get 
$L_n=\left \lceil (n+1-c-2d)/2c 
\right \rceil -1$ and $M_n=\left \lfloor (n-d-1)/c \right \rfloor -1$. Notice 
that $$f_w(n+2c)-f_w(n)=\sum_{i=1}^{2c} M_{n+i}- \sum_{i=1}^{2c} L_{n+i-1}+2c=$$ 
$$\left \lfloor (n-d)/c \right \rfloor +...+
\left \lfloor (n+2c-d-1)/c \right \rfloor -\left \lceil (n+1-c-2d)/2c 
\right \rceil -...-$$ $$\left \lceil (n+c-2d)/2c \right \rceil +2c = (n-d)+(n-d+c)-
(n+c-2d) +2c=n+2c.$$ So $f_w(n+2c)=f_w(n)+n+2c.$ 
\end{proof}

\section{Necessary and sufficient conditions}
\label{sec:cond}

In this section we give necessary and sufficient conditions for a function 
$f:\Z _+ \rightarrow \Z _+$ to be the subword complexity of an infinite gap 
increasing word. We use the notations $\Delta f (n) = f(n+1)-f(n)$ and $\Delta ^2 f (n) = \Delta f(n+1)- \Delta f(n)$.

The next proposition gives a necessary condition for a function to be
the subword complexity function of an infinite binary word whose gap 
function in injective. A stronger condition is obtained for gap 
increasing words.  

\begin{proposition} Let $w$ be an infinite binary word with gap function 
$g$  and subword complexity function $f$. If $g$ is injective, then 
$1 \leq \Delta f(n) \leq n+1.$ If $g$ is strictly increasing, then 
$1 \leq \Delta f(n) \leq \left \lfloor n/2 \right \rfloor + 1.$
\label{pr:simnec}
\end{proposition}

\begin{proof} Clearly $\Delta f(n)$ is the number of right special factors
of length $n$ of $w$.

Let $g$ be injective, then it follows from the proof of
Proposition~\ref{pr:inj} that the number of right special factors of $w$ of
length $n$ equals  one plus the number of subwords of $w$ of form $v=0^{x}10^{y}1$,
for some $x \geq 0$ and $y \geq 0$ such that $x+y=n-1$. Thus $1 \leq \Delta f(n) \leq n+1.$

If $g$ is strictly increasing, then $0 \leq x <y$, hence $x$ can 
take at most $\left \lfloor n/2 \right \rfloor$ different values and 
$1 \leq \Delta f(n) \leq \left \lfloor n/2 \right \rfloor + 1.$
\end{proof}

\begin{corollary} Let $w$ be an infinite gap increasing word with subword
complexity function $f$, then  $n+1 \leq f(n) \leq  \left \lceil n/2
\right \rceil \left \lfloor n/2 \right \rfloor +\left \lceil n/2
\right \rceil +1.$ 
\label{cr:bds} 
\end{corollary}

\begin{proof} This follows from $f(n)=2+\sum_{i=1}^{n-1} \Delta f(i)$ and  
$1 \leq \Delta f(i) \leq \left \lfloor i/2 \right \rfloor + 1$ for all $i$.
\end{proof}

\begin{remark} The lower bound in Corollary~\ref{cr:bds} is exact: for 
any positive integer $N$ there is an infinite gap
increasing word $w$ with the property that $f(n+1)=n+2$ for $0 \leq n \leq
N$ ($w$ could be any gap increasing word with prefix $0^k1$ where $k \geq
N-1$). However there is no infinite gap increasing word such that
$f(n+1)=n+2$ for all $n$, that is an infinite gap increasing word cannot be
Sturmian. The upper bound in Corollary~\ref{cr:bds} is
also exact (this will be shown in Lemma~\ref{lm:exbd}).
\end{remark}

\begin{definition} We will say that an infinite gap increasing word $w$ has a double gap of length $n$ if there exists an integer $m \geq 0$ such that either $10^{m}10^{n-m}1$ is a subword of $w$ or $0^{m}10^{n-m}1$ is a prefix of $w$, in other words there exist two consecutive "gaps" (blocks of $0$'s between $1$'s or before the first $1$) whose length sum is $n$.   
\end{definition}

\begin{lemma} Let $f$ be a subword complexity function of an infinite 
gap increasing word. Then, for each $n \in \Z _+$, one of the following statements is true:

a) $10^n1$ is a subword of $w$ and $w$ does not have a double gap of length $n-1$, then
${\Delta }^2 f(n) = 1$,

b) $10^n1$ is not a subword of $w$ and $w$ has a double gap of length $n-1$, then ${\Delta }^2 f(n) = -1$, 

c) $10^n1$ is a subword of $w$ and $w$ has a double gap of length $n-1$, then ${\Delta }^2 f(n) = 0$,

d) $10^n1$ is not a subword of $w$ and $w$ does not have a double gap of length $n-1$, then ${\Delta }^2 f(n) = 0$.
\label{lm:cond}
\end{lemma}

\begin{proof} Clearly $\Delta f(n)=f(n+1)-f(n)$ is the number of right special factors of $w$ of length $n$. By the proof of Proposition~\ref{pr:simnec}, the number of right special factors of $w$ of length $n$ equals the number of subwords of $w$ of form $0^{x}10^{n-x-1}$, where $x \geq 0$, plus one. Thus ${\Delta }^2 f(n) = \Delta f(n+1)-\Delta f(n)$ is the difference between the number of subwords of $w$ of form $0^{x}10^{n-x}1$ and the number of subwords of $w$ of form $0^{x}10^{n-x-1}1$.     
 
For any $n \in \Z _+$, let $C(n)$ denote the set of subwords of $w$ of form
$0^{x}10^{n-x-1}1$. Then ${\Delta }^2 f(n) = |C(n+1)|-|C(n)|$.

For every $v \in C(n+1)$, let $\phi (v)$ be the suffix of length
$n+1$ of v. If $x \not= 0$, then $\phi (0^x10^{n-x}1)) = 0^{x-1}10^{n-x}1 \in C(n)$. If $x=0$, then $\phi (10^n1) = 0^n1 \notin C(n)$. We conclude that $\phi$ maps all elements of $C(n+1)$ but $10^n1$ (if it happens to be a subword of $w$) to $C(n)$. 

For every element $u = 0^{x}10^{n-x-1}1 \in C(n)$, $\phi ^{-1} (u)$ exists 
if and only if  $0^{x+1}10^{n-x-1}1$ is a subword of $w$. So $\phi ^{-1} (u)$ does not exist if either $u=10^x10^{n-x-1}1$ is a subword of $w$ or $u=0^{x}10^{n-x-1}1$ is a prefix of $w$, which happens when $w$ has a double gap of length $n-1$.

In case (a), when $10^n1$ is a subword of $w$ and $w$ does not have a double gap of length $n-1$, $\phi ^{-1}:C(n) \rightarrow C(n+1)$ is a well defined function and 
$C(n+1)$ has one more element than $C(n)$, so ${\Delta }^2 f(n) = 1$.

In case (b) $\phi :C(n+1) \rightarrow C(n)$ is an injective function and $C(n)$ 
has one more element than $C(n+1)$, so ${\Delta }^2 f(n) = -1$.
 
In case (c) $\phi$ is a bijection, hence ${\Delta }^2 f(n) = 0$.

In case (d) $\phi$ maps all but one element of $C(n+1)$ to $C(n)$ and 
$\phi ^{-1}$ exists for all but one element of $C(n)$. Therefore the cardinalities of $C(n+1)$ and $C(n)$ are the same and ${\Delta }^2 f(n) = 0$.  
\end{proof}

\begin{corollary} [Necessary condition for a function to be the subword
complexity function of an infinite gap increasing word] If $f$ is a subword
complexity function of an infinite gap increasing word, then $f(1)=2$, $f(2)=3$ and $| { \Delta }^2 f(n) | \leq 1$.
\label{cor:simnec}
\end{corollary}
\begin{proof} It follows from Lemma~\ref{lm:cond} that $| { \Delta }^2 f(n) | \leq 1$.
$f(2)=3$ because $00$, $01$ and $10$ are necessarily subwords of an infinite gap increasing word, while $11$ cannot be a subword of an infinite gap increasing word.  
Indeed, if an infinite binary words starts with $11$, then the gap function g of this word satisfies $g(1)=g(2)=1$ and thus $g$ is not strictly increasing. 
\end{proof}

Next we introduce some new terminology that will enable us to develop a method to find out if a given function is a subword complexity function of an infinite gap increasing word. 

It will be convenient to think of sets of positive integers as of strictly increasing sequences of positive integers. These sequences can be finite, infinite or even \textbf {empty}.
We will make the following abuse of notation: if $\{ x_i  \} _{i \geq 1}$ and
$\{ y_i  \} _{i \geq 1}$ are two strictly increasing sequences that have no elements
in common, then $\{ x_i  \} _{i \geq 1} \sqcup \{ y_i  \} _{i \geq 1}$ is the strictly increasing sequence that consists of all the elements of both sequences. 
If $\{ y_i  \} _{i \geq 1}$ is empty, then $\{ x_i  \} _{i \geq 1} \sqcup \{ y_i  \} _{i \geq 1} = \{ x_i  \} _{i \geq 1}$. 

\begin{example} $$\{ 2i  \} _{i=1}^{\infty} \sqcup \{ 2i-1  \} _{i=1}^4= \{ 1,2,3,4,5,6,7,8,10,12,14,16,18,...\}= \{ i \} _{i=1}^8 \sqcup \{ 2i  \} _{i=5}^{\infty} .$$ 
\end{example}

\begin{definition}  Let $\{ x_i  \} _{i \geq 1}$ be a strictly increasing sequence of positive integers and $p<x_1$ be a non-negative integer. Then $\sigma_p ( \{ x_i  \} _{i \geq 1} )$ will denote the strictly increasing sequence $\{ p+x_1+1,x_1+x_2+1,x_2+x_3+1,... \} $.
\end{definition}

\begin{example} $$\sigma_1 ( \{ 2i  \} _{i=1}^{\infty} ) = \{ 4,7,11,15,... \} = \{ 4 \}
\sqcup \{ 4i+3  \} _{i=1}^{\infty} .$$
\end{example}

\begin{theorem} [Necessary and sufficient condition for a function to be the
subword complexity function of an infinite gap increasing word] Let $f: \Z _+
\rightarrow \Z _+$ satisfy the necessary condition of Proposition~\ref{pr:simnec},
that is $f(1)=2$, $f(2)=3$ and $| { \Delta }^2 f(n) | \leq 1$. Set $\{ a_i  \} _{i \geq 1}$
to the sequence obtained by arranging the elements of the set 
$\{ n \in \Z _+ |\Delta ^2 f(n) = 1 \}$ in strictly increasing order. 
Similarly $\{ b_i  \} _{i \geq 1}$ is obtained by arranging the elements of $\{ n \in \Z _+ | \Delta ^2 f(n) = -1 \}$
in strictly increasing order. Then $f$ is the subword complexity function of an infinite gap
increasing word if and only if $\{ a_i  \} _{i \geq 1}$ is not empty and there exist an integer $p$, 
$0 \leq p < a_1$, and a strictly increasing sequence of positive integers $\{ c_i  \} _{i \geq 1}$, 
with no elements in common with either $\{ a_i  \} _{i \geq 1}$ or $\{ b_i  \} _{i \geq 1}$, 
such that $\{ a_i  \} _{i \geq 1} \sqcup \{ c_i  \} _{i \geq 1}$ is infinite and 
\begin{equation}
\sigma_p (\{ a_i  \} _{i \geq 1} \sqcup \{ c_i  \} _{i \geq 1}) = 
\{ b_i  \} _{i \geq 1} \sqcup \{ c_i  \} _{i \geq 1}.
\label{eq:seqset}
\end{equation} 
\label{th:suf}
\end{theorem} 

\begin{proof} 
($\rightarrow$) Let $f$ be the subword complexity function of the infinite gap
increasing word $w = 0^q10^{j_1}10^{j_2}10^{j_3}1...$, where $\{ j_i  \} _{i=1}^{\infty}$ is a strictly increasing sequence of positive integers and
the integer $q$ is in the range $0 \leq q < j_1$. We want to show that $\{ a_i  \} _{i \geq 1}$ is not empty and that there exist
an integer $p$, $0 \leq p < a_1$, and a strictly increasing sequence of positive integers $\{ c_i  \} _{i \geq 1}$, with no elements in common with either $\{ a_i  \} _{i \geq 1}$ or $\{ b_i  \} _{i \geq 1}$, such that $\{ a_i  \} _{i \geq 1} \sqcup \{ c_i  \} _{i \geq 1}$ is infinite and eq.~(\ref{eq:seqset}) is satisfied.

Notice that $10^n1$ is a subword of $w$ if and only if $n$ is an element of $\{ j_i  \} _{i=1}^{\infty}$. Also $w$ has a double gap of length $n-1$ if and only if $n-1=q+ j_1$ or $n-1=j_i+ j_{i+1}$ for some $i$, which is equivalent to saying that $n$ is an element of $\sigma _q(\{ j_i  \} _{i=1}^{\infty})$. Let $\{ c_i  \} _{i \geq 1}$ be the strictly increasing sequence of integers $n$ such that $10^n1$ is a subword of $w$ and $w$ has a double gap of length $n-1$. Clearly $\{ c_i  \} _{i \geq 1}$ has no elements in common with either $\{ a_i  \} _{i \geq 1}$ or $\{ b_i  \} _{i \geq 1}$. It also follows from Lemma~\ref{lm:cond} that $\{ j_i  \} _{i=1}^{\infty}=\{ a_i  \} _{i \geq 1} \sqcup \{ c_i  \} _{i \geq 1}$ and $\sigma_q (\{ j_i  \} _{i=1}^{\infty}) = \{ b_i  \} _{i \geq 1} \sqcup \{ c_i  \} _{i \geq 1}$. Set $p=q$. We proved that $\{ a_i  \} _{i \geq 1} \sqcup \{ c_i  \} _{i \geq 1}$ is infinite and eq.~(\ref{eq:seqset}) is satisfied. We still need to prove that $\{ a_i  \} _{i \geq 1}$ is not empty and $0 \leq p < a_1$. To prove both statements it is enough to prove that $j_1=a_1$. Suppose $j_1 \ne a_1$, then $j_1 =c_1$, which is impossible because $\{ c_i  \} _{i \geq 1}$ is a subsequence of $\sigma_q (\{ j_i  \} _{i=1}^{\infty})$ and the first (and the least) element of $\sigma_q (\{ j_i  \} _{i=1}^{\infty}$ is $p+j_1+1>j_1$.  

($\leftarrow$) Let the function $f$ satisfy the necessary condition of Proposition~\ref{pr:simnec} and the sequences $\{ a_i  \} _{i \geq 1}$ and $\{ b_i  \} _{i \geq 1}$ be as described in the hypothesis. Further suppose that $\{ a_i  \} _{i \geq 1}$ is not empty and there exist an integer $p$, $0 \leq p < a_1$, and a strictly increasing sequence of positive integers $\{ c_i  \} _{i \geq 1}$, 
with no elements in common with either $\{ a_i  \} _{i \geq 1}$ or $\{ b_i  \} _{i \geq 1}$, such that $\{ a_i  \} _{i \geq 1} \sqcup \{ c_i  \} _{i \geq 1}$ is infinite and 
eq.~(\ref{eq:seqset}) is satisfied. We want to show that $f$ is the subword complexity function of a gap increasing word.

Denote the sequence $\{ a_i  \} _{i \geq 1} \sqcup \{ b_i  \} _{i \geq 1}$ by $\{ j_i  \} _{i=1}^{\infty}$. It follows from eq.~(\ref{eq:seqset}) that $j_1=a_1$. Thus $p<j_1$ and the infinite word $v = 0^p10^{j_1}10^{j_2}10^{j_3}1...$ is gap increasing. Let $f_v$ be the subword complexity of $v$. We will show that $f_v=f$. 

Since $v$ is gap increasing, $f_v(1)=2=f(1)$ and $f_v(2)=3=f(2)$. It is therefore enough to show that ${ \Delta }^2 f_v(n)={ \Delta }^2 f(n)$ for all positive integers $n$. It follows from eq.~(\ref{eq:seqset}) that all elements that $\{ j_i  \} _{i=1}^{\infty}$ and $\sigma_q (\{ j_i  \} _{i=1}^{\infty})$ have in common form the sequence $\{ c_i  \} _{i \geq 1}$. By Lemma~\ref{lm:cond} ${ \Delta }^2 f_v(a_i)=1$
for all the elements of $\{ a_i  \} _{i \geq 1}$, ${ \Delta }^2 f_v(b_i)=-1$ for all the elements of $\{ b_i  \} _{i \geq 1}$, ${ \Delta }^2 f_v(n)=0$ for all positive integers $n$ that are neither in $\{ a_i  \} _{i \geq 1}$ nor in $\{ b_i  \} _{i \geq 1}$. This means that ${ \Delta }^2 f_v(n)={ \Delta }^2 f(n)$ for all positive integers $n$.
\end{proof}

Based on Theorem~\ref{th:suf} we give a practical method of determining if a
function $f: \Z _+ \rightarrow \Z _+$ is the subword complexity function of an infinite gap increasing word.

From now on let's fix a function $f: \Z _+ \rightarrow \Z _+$ such that $f(1)=2$, $f(2)=3$ and $| { \Delta }^2 f(n) | \leq 1$. Further let's assume that the sequences $\{ a_i  \} _{i \geq 1}$ and $\{ b_i  \} _{i \geq 1}$, defined in Theorem~\ref{th:suf}, are already computed and $\{ a_i  \} _{i \geq 1}$ is not empty. We need to find a way to determine if there exists an integer $p$, $0 \leq p < a_1$ and a strictly increasing sequence of positive integers $\{ c_i  \} _{i \geq 1}$, with no elements in common with either $\{ a_i  \} _{i \geq 1}$ or $\{ b_i  \} _{i \geq 1}$, such that eq.~(\ref{eq:seqset}) holds. The method involves testing every $p$, $0 \leq p < a_1$.

\textbf{Fix some integer} $\bf{p}$, $\bf{0 \leq p < a_1}$. There exists a sequence $\{ c_i  \} _{i \geq 1}$ that satisfies the description above if and only if there exists a sequence $\{ j_i  \} _{i=1}^{\infty} = \{ a_i  \} _{i \geq 1} \sqcup \{ c_i  \} _{i \geq 1}$, which satisfies the following three properties:

P1) $j_1=a_1$ and $\{ a_i  \} _{i \geq 1}$ is a subsequence of $\{ j_i  \} _{i=1}^{\infty}$;

P2) $\{ b_i  \} _{i \geq 1}$ is a subsequence of $\{ p+j_1+1,j_1+j_2+1,j_2+j_3+1,... \}$;

P3) Every element of the sequence $\{ p+j_1+1,j_1+j_2+1,j_2+j_3+1,... \}$, that is not an element of the sequence $\{ b_i  \} _{i \geq 1}$, is an element of the sequence $\{ j_i  \} _{i=1}^{\infty}$. No element of the sequence $\{ p+j_1+1,j_1+j_2+1,j_2+j_3+1,... \}$ is an element of $\{ a_i  \} _{i \geq 1}$.

It follows from the proof of the Theorem~\ref{th:suf} that if the sequence $\{ j_i  \} _{i=1}^{\infty}$ with the properties P1-P3 above exists, then the word $w = 0^p10^{j_1}10^{j_2}10^{j_3}1...$ has the subword complexity function $f$.

\begin{proposition} Properties P1-P3 imply the following recursive construction of $\{ j_i  \} _{i=1}^{\infty}$. 

R1) The initial conditions:

$j_1=a_1$ and $\{ a_i  \} _{i \geq 1}$ is a subsequence of $\{ j_i  \} _{i=1}^{\infty}$;

R2) The recurrence:

Every element of $\{ j_i  \} _{i=1}^{\infty}$ that falls in the interval 
$(a_s,a_{s+1}) \cap (b_t,b_{t+1})$ (if $a_s$ is the last element of $\{ a_i  \} _{i \geq 1}$, then $(a_s,a_{s+1})=(a_s,\infty)$; if $b_t$ is the last element of $\{ b_i  \} _{i \geq 1}$, then $(b_t,b_{t+1})=(b_t,\infty)$) is computed by:

$$j_i=j_{i-s+t} + j_{i-s+t-1}+1.$$  
 
T) The test: The sequence $\{ j_i  \} _{i=1}^{\infty}$ that is recursively computed using R1-R2 satisfies properties P1-P3 if and only if $\{ b_i  \} _{i \geq 1}$ is a subsequence of $\sigma _q(\{ j_i  \} _{i=1}^{\infty}) = \{ p+j_1+1,j_1+j_2+1,j_2+j_3+1,... \}$ and $\{ a_i  \} _{i \geq 1}$ has no elements in common with $\sigma _q(\{ j_i  \} _{i=1}^{\infty})$.
\label{pr:C}
\end{proposition}

\begin{corollary} Let $f$ be the subword complexity function of an infinite gap increasing word $w$ with prefix $0^p1$ for some nonnegative integer $p$. Then $w$ is the only gap increasing word with prefix $0^p1$ and subword complexity function $f$. 
\label{cor:uniC}
\end{corollary}

\begin{proof} For a fixed function $f$ (and therefore fixed $\{ a_i  \} _{i \geq 1}$ and $\{ b_i  \} _{i \geq 1}$) and a fixed $p$, the recursive construction given in Proposition~\ref{pr:C} gives a unique sequence $\{ j_i  \} _{i=1}^{\infty}$. Thus there is a unique gap increasing word $0^p10^{j_1}10^{j_2}10^{j_3}1...$ with the subword complexity function $f$.
\end{proof}

\begin{remark} By Corollary~\ref{cor:uniC}, for any $f: \Z _+
\rightarrow \Z _+$ and any nonnegative integer $p$ there is at most one
infinite gap increasing word with prefix $0^p1$ with subword complexity
function $f$. However it is not true in general that a function $f: \Z _+
\rightarrow \Z _+$ cannot be the subword complexity function of two distinct
infinite gap increasing words. 

Consider the following two gap increasing words:
$$u=010^310^510^910^{15}10^{25}...$$ with gap function $g_u$
defined by $g_u(1)=2$, $g_u(2)=4$ and $g_u(i)=g_u(i-2)+g_u(i-1)$ for $i \geq
3$,  and
$$v=0^210^310^610^{10}10^{17}10^{28}...$$ with gap function $g_v$
defined by $g_v(1)=3$, $g_v(2)=4$ and $g_v(i)=g_v(i-2)+g_u(i-1)$ for $i \geq
3$.

To prove that $u$ and $v$ have the same subword complexity it suffices to prove that the pair of sequences $\{ a_i  \} _{i \geq 1}$ and $\{ b_i  \} _{i \geq 1}$, that correspond to $u$, and the pair of sequences  $\{ a'_i  \} _{i \geq 1}$ and $\{ b'_i  \} _{i \geq 1}$, that correspond to $v$, are the same.

By Lemma~\ref{lm:cond} the sequences $\{ a_i  \} _{i \geq 1}$ and $\{ a'_i  \} _{i \geq 1}$ consist of only one element $a_1=a'_1=3$; while the sequences $\{ b_i  \} _{i \geq 1}$ and $\{ b'_i  \} _{i \geq 1}$ are empty.

Hence $u$ and $v$ have the same subword complexity function $f(n)=2+\sum_{i=1}^{n-1} \Delta f(n) $, where $\Delta f(i)=1$ for $i \leq 3$ and $\Delta f(i)=2$ for $i \geq 4$. Thus $f(n)=n+1$ for $n \leq 4$ and $f(n)=2n-3$ for $n \geq 5$.   
\label{rm:uv}
\end{remark}

The next two examples illustrate how one can use Proposition~\ref{pr:C} to determine if a function $f:\Z _+ \rightarrow \Z _+$ is the subword complexity function of an infinite 
gap increasing word. 

\begin{example} We want to find out if there exists an integer $q>1$ for which $$f(n)=  
(\left \lfloor \frac {n}{q} \right \rfloor +1)(n - \frac {q}{2} 
\left \lfloor \frac {n}{q} \right \rfloor) +1$$ is a subword complexity function of an infinite gap increasing word. 

First we check that the necessary condition of Proposition~\ref{pr:simnec} are satisfied. Indeed $f(1)=2$, $f(2)=3$ and $$\Delta f(n)=f(n+1)-f(n)=\left \lceil \frac {n+1}{q} \right \rceil.$$
Next we find the sequences $\{ a_i  \} _{i \geq 1}$ and $\{ b_i  \} _{i \geq 1}$
defined in Theorem~\ref{th:suf}. The sequence $\{ b_i  \} _{i \geq 1}$ is empty.
The sequence $\{ a_i  \} _{i=1}^{\infty}$ is given by $a_i=qi-1$.

Next we compute the sequence $\{ j_i  \} _{i=1}^{\infty}$
using the recurrence in Proposition~\ref{pr:C}. For any $p$, $0 \leq p \leq q-2$, $\{ j_i  \} _{i=1}^{\infty} = \{ q-1,p+q,2q-1,p+2q,3q-1,p+3q,4q-1,p+4q,5q-1,p+5q,... \} $. By Proposition~\ref{pr:C}
the function $f(n)$ is the subword complexity function of a gap increasing word because $\{ j_i  \} _{i=1}^{\infty}$ satisfies test T, that is $\{ a_i  \} _{i=1}^{\infty}=\{ qi-1  \} _{i=1}^{\infty}$ has no elements in common with $\sigma _q(\{ j_i  \} _{i=1}^{\infty})=
\{ p+q,p+2q,p+3q,p+4q,p+5q,... \}$ for any $p$, $0 \leq p \leq q-2$.

This shows that for any integer $q>1$ the function $$f(n)=  
(\left \lfloor \frac {n}{q} \right \rfloor +1)(n - \frac {q}{2} 
\left \lfloor \frac {n}{q} \right \rfloor) +1$$ is the subword complexity function of exactly $q-1$ distinct infinite gap increasing words.
\label{ex:yes}   
\end{example}

\begin{example} We want to find out if there exists an integer $q>2$ for which $$f(n)= n - \frac {q}{2} 
\left \lfloor \frac {n-1}{q} \right \rfloor ^2 + (n-1 - \frac {q}{2}) 
\left \lfloor \frac {n-1}{q} \right \rfloor +1$$ is a subword complexity function of an infinite gap increasing word. 

One can check that $f(1)=2$, $f(2)=3$ and 
$$\Delta f(n)=f(n+1)-f(n)=\left \lfloor \frac {n-1}{q} \right \rfloor+1.$$
Next we find the sequences $\{ a_i  \} _{i \geq 1}$ and $\{ b_i  \} _{i \geq 1}$
defined in Theorem~\ref{th:suf}. The sequence $\{ b_i  \} _{i \geq 1}$ is empty.
The sequence $\{ a_i  \} _{i=1}^{\infty}$ is given by $a_i=qi$.

We intend to check if the sequence $\{ j_i  \} _{i=1}^{\infty}$ that is computed using R1-R2 in Proposition~\ref{pr:C} also satisfies the test T for at least one $p$, $0 \leq p \leq q-1$. For any $p$, $\{ j_i  \} _{i=1}^{\infty} = \{ q,p+q+1,2q,p+2q+2,3q,p+3q+2,4q,p+4q+3,5q,p+5q+3,... \} $. To satisfy test T the sequence $\{ a_i  \} _{i=1}^{\infty}=\{ qi  \} _{i=1}^{\infty}$ should not have elements in common with $\sigma _q(\{ j_i  \} _{i=1}^{\infty}) = \{ p+iq+r_i \} _{i=1}^{\infty}$, where $\{ r_i  \} _{i=1}^{\infty}$ is the non-decreasing sequence in which every positive integer $n$ occurs $2^{n-1}$ times, that is $\{ r_i  \} _{i=1}^{\infty} = \{ 1,2,2,3,3,3,3,4,4,4,4,4,4,4,4,5,... \}$. For $k$ large enough (more precisely $k=2^{q-p-1}$) $r_k=q-p$ and thus the $k^{th}$ element of $\sigma _q(\{ j_i  \} _{i=1}^{\infty})$ is $p+kq+q-p=(k+1)q=a_{k+1}$. Thus the test is not satisfied for any $p$, $0 \leq p \leq q-1$.

This shows that there is no integer $q>2$ for which $$f(n)= n - \frac {q}{2} 
\left \lfloor \frac {n-1}{q} \right \rfloor ^2 + (n-1 - \frac {q}{2}) 
\left \lfloor \frac {n-1}{q} \right \rfloor +1$$ is a subword complexity function of an infinite gap increasing word. 
\label{ex:no}  
\end{example}

\begin{remark} By Remark~\ref{rm:uv} and Example~\ref{ex:yes} there can exist several different infinite gap increasing words with the same subword complexity function f. However we conjecture that, if there exists an $n$ such that $\Delta^2 f(n)=-1$, then there exists at most one infinite gap increasing word with subword complexity function $f$.  
\end{remark}

\section{A geometric representation of the subword complexity of gap increasing words}
\label{sec:aff}
 
In this section we find a geometric representation of the subword complexity function of a gap increasing word, which allows us to compute for each positive integer $n$ the exact upper bound of the subword complexity function $f_w(n)$ over all infinite gap increasing words $w$. It will turn out that for all $n$ the value of $f_w(n)$ is maximized by the same $w$. 
  
We say that an infinite word $w$ has subword complexity $\phi (n)$ {\it ultimately} if there exists a positive integer $N$, such that for all integers $n \geq N$, the subword complexity function of $w$ is $\phi (n)$. In this section we determine for which integers $a$ and $b$ there is an infinite gap increasing word $w$ with the subword complexity function $an+b$ ultimately. 

\begin{theorem} Let $w=0^{n_0}10^{n_1}10^{n_2}10^{n_3}...$ be an infinite gap increasing word. Let $f_w$ be the subword complexity function of $w$. For a fixed  $n>n_1$, let the partition $\nu _n$ be the partition whose all parts are all $n_i+1$, such that $n_i+1\leq n$. In drawing the diagram of partition $\nu _n$ we adopt the French convention, that is the bottom row is the longest row and the left-most column is the longest column. Let $r(n)+1$ be the number of parts of $\nu _n$. The diagram of $\nu _n$ is contained in the diagram of $(n+1)^{r(n)+1}$, which is an $(n+1) \times (r(n)+1)$ rectangle (as shown in fig.~\ref{irinafig}).

Consider the boundary line between $\nu _n$ and its complement in $((n+1)^{r(n)+1})$ (marked thickly in fig.~\ref{irinafig}). Index the rows of the diagram of $((n+1)^{r(n)+1})$ from $0$ (top) to $r(n)$ (bottom). Let $l(n)$ be the maximum row index (if it exists) for which the portion of the $(l(n)-1)^{th}$ row to the left of the boundary line is less than the portion of the $l(n)^{th}$ row to the right of the boundary line. If such row does not exist set $l(n)=0$. If $l(n)>0$, shade the portion of the diagram of $\nu _n$ above the $l(n)^{th}$ row and, if $l(n)<r(n)$, shade the portion of the complement of $\nu _n$ in $((n+1)^{r(n)+1})$ below the $l(n)^{th}$ row. Then $f_w(n+1)=n+2+{\text {shaded area}}$. 

\begin{figure}[!ht]
\begin{center}
\includegraphics[scale=0.9]{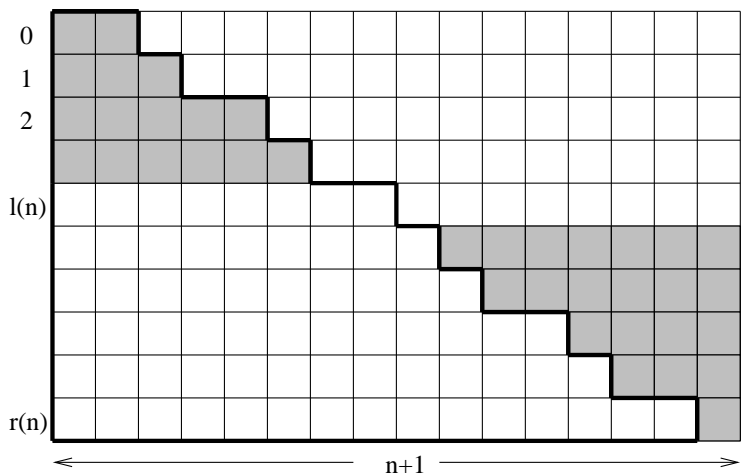}
\caption{The diagram of $\nu _n$ (bounded by the thick line) in the 
$(n+1) \times (r(n)+1)$ rectangle}
\label{irinafig}
\end{center}
\end{figure}

\label{th:figure}

Remark: for $n \leq n_1$ the subword complexity function of $w$ is $f_w(n+1)=n+2$.
\end{theorem}

\begin{proof}
It follows from Proposition~\ref{pr:inj} that the number of right special factors of $w$ of length $k>n_0$ is $s_w(k)=1 + r(k) - l(k)$, where $r(k)$ is the number of $i \geq 1$ such that $n_i \leq k-1$ and $l(k)$ is the number of $i \geq 1$ such that $n_{i-1}+n_i \leq k-2$ (it will be shown later that $l(k)$ can be defined the way it was defined in the hypothesis of Theorem~\ref{th:figure}). The $(n+1)$th subword complexity of $w$ is
\begin{equation*}
f_w(n+1)=f_w(1)+ \sum _{k=1}^n s_w(k)= 2+ \sum _{k=1}^n (1 + r(k) - l(k))=
\end{equation*} 
\begin{equation*}
n+2+ \sum _{i=1}^{r(n)} (n-n_i)- \sum _{i=1}^{l(n)} (n-n_{i-1}-n_i-1)=
\end{equation*}
\begin{equation*}
n+2+\sum _{i=1}^{l(n)}(n_{i-1}+1) + \sum _{i=l(n)+1}^{r(n)}(n-n_i).
\end{equation*}

Thus
\begin{equation}
f_w(n+1)=n+2+\sum _{i=0}^{l(n)-1}(n_i+1) + \sum _{i=l(n)+1}^{r(n)}(n-n_i).
\label{eq:diag}
\end{equation}

Consider the partition $\nu _n = (n_0+1,n_1+1,...,n_{r(n)}+1)$, where all parts $n_i+1 \leq n$. By eq.~\ref{eq:diag} $f_w(n+1)$ equals $n+2$ plus the area of the figure that consists of the portion of $\nu _n$ above the $l(n)$th row and the portion of the complement of $\nu _n$ in $((n+1)^{r(n)+1})$ below the $l(n)$th row (shaded area in fig.~\ref{irinafig}). Clearly $l(n)=0$ if $n_0+1 \geq n-n_1$. Otherwise $l(n)$ is maximal with property $n_{l(n)-1}+1 < n-n_{l(n)}$. In the last case we can define $l(n)$ in terms of the diagram as the maximum integer for which the part of $\nu _n$ in the $(l(n)-1)^{th}$ row is less than the part of the complement of $\nu _n$ in the $l(n)^{th}$ row. 
\end{proof}

\begin{lemma} The exact upper bound on the subword complexity function $f(n)$ of an 
infinite gap increasing word is $\left \lceil n/2 \right \rceil \left \lfloor n/2 
\right \rfloor + \left \lceil n/2 \right \rceil +1$.
\label{lm:exbd}
\end{lemma}

\begin{proof} Clearly the shaded area in fig.~\ref{irinafig} is maximized by
the gap increasing word $w$ with $n_i=i$ for $i \geq 0$. The gap function of 
$w$ is $g(n)=n$ and it follows from Lemma~\ref{lm:cn+d} that f(n+2)=f(n)+n+2.
Solving this recurrence with initial conditions f(1)=2 and f(2)=3, we get  
$$f(n)= \left \lceil n/2 \right \rceil \left \lfloor n/2 \right \rfloor + \left 
\lceil n/2 \right \rceil +1.$$
\end{proof}

The diagram of partition $\nu _n$ in $((n+1)^{r(n)+1})$ (fig.~\ref{irinafig}), used to compute $f_w(n+1)$, is a representation of a prefix of $w$ of length equal to $|\nu _n|$ (the weight of $\nu _n$). We will call such a diagram a gap increasing prefix diagram. 

A diagram of any partition $\mu$ with distinct parts, in a rectangle with the vertical side of length equal to the number of parts in $\mu$ and with the horizontal side of length greater then the largest part of $\mu$, can be thought as a gap increasing prefix diagram. Let $k$ be the length of the horizontal side of the rectangle. Theorem~\ref{th:figure} gives a method for computing the $k^{th}$ value of the subword complexity function of some infinite gap increasing word whose prefix is represented by $\mu$ and whose suffix (one that follows after the prefix represented by $\mu$) does not contain subwords $10^n1$ for $n<k-1$. We will call this value the complexity $\bar {f}(k)$ of the gap increasing prefix diagram.

\begin{theorem} Let $a$ and $b$ be two integers. There exists a gap increasing word $w$ with the subword complexity function $f_w(n)=an+b$ ultimately if and only if $a \geq 2$.
\label{th:anb} 
\end{theorem}
 
\begin{proof} ( $\rightarrow$ ) Let $a \geq 2$ and $b$ be two integers. To show the existence of an infinite gap increasing word $w$ with the subword complexity $an+b$ ultimately, we first build a prefix of $w$ and then show the recursive construction of the infinite suffix of $w$.   

Our first goal is to build a prefix $v$ of $w$ with the following property: there exists a positive integer $K$, such that any word $vu$, where $u$ is an infinite suffix that makes $vu$ gap increasing and $u$ has no subwords $10^n1$ with $n<K-1$, has $K$th subword complexity $f_{uv}(K)=aK + b$ and the number of right special factors of $vu$ of length $K-1$ is $s_{uv}(K-1)=a$. To build $v$, it is enough to construct its gap increasing diagram $D$ with complexity $\bar {f}_D (K)=aK+b$ such that the diagram obtained by deleting the right-most column of $D$ has complexity $\bar {f}_D (K) -a$.   

\begin{figure}
\begin{center}
\includegraphics[scale=0.8]{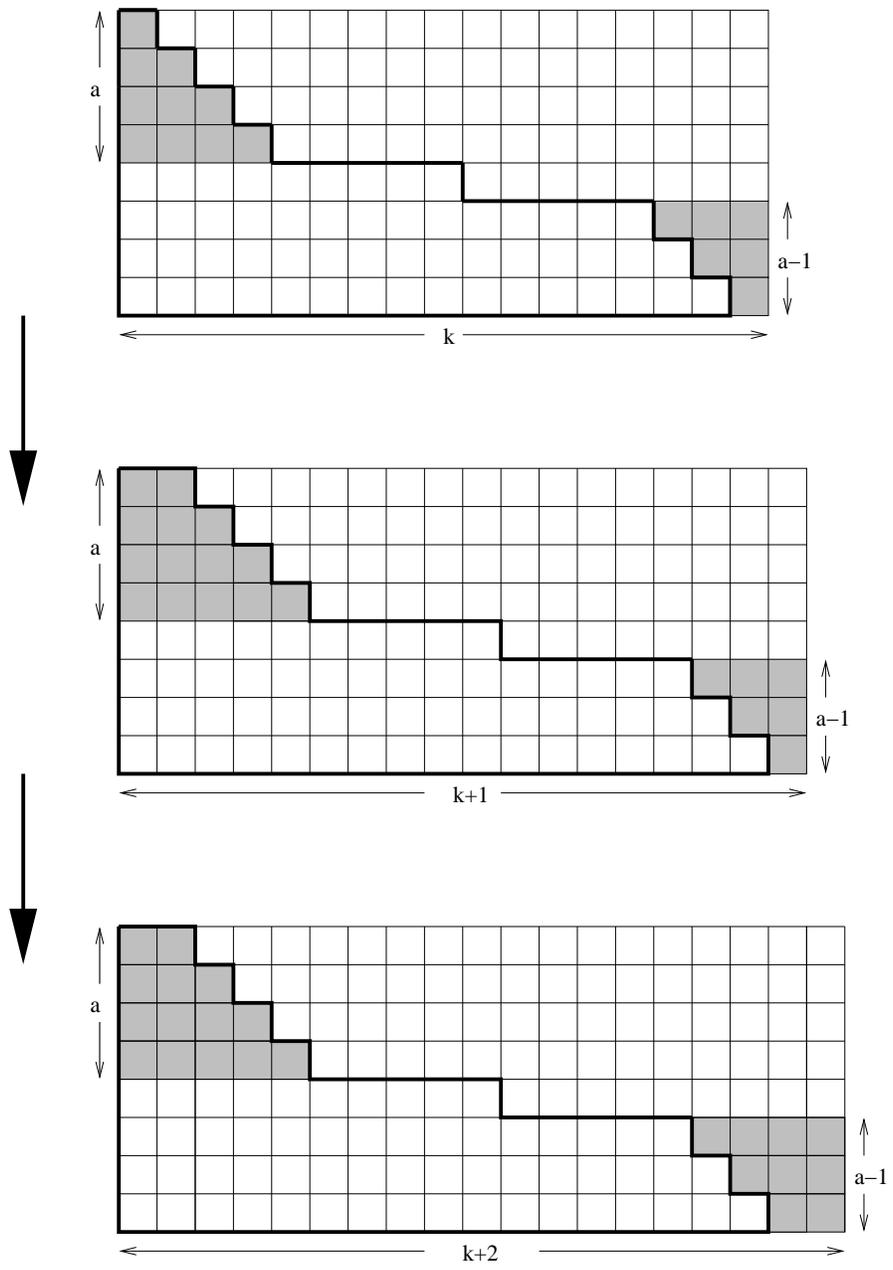}
\caption{The steps of the algorithm for constructing the diagram of a prefix of $w$}
\label{omarfig}
\end{center}
\end{figure}

Since $a>1$, there exists an integer $k_0 \geq 2a+4$ such that $k_0+a^2 < ak_0+b$. The diagram of partition $(1,2,3,...,a-1,a,\left \lfloor k_0/2 \right 
\rfloor , k_0-a+1,k_0-a+2,k_0-a+3,...,k_0-1)$ with $2a$ parts, in a rectangle with the vertical side of length $2a$ and the horizontal side of length $k_0$, is a gap increasing partition diagram, which we denote by $D_0$. By Theorem~\ref{th:figure}, the  ${k_0}^{th}$ complexity $\bar {f_0} (k_0)$ of this digram  is $k_0+1$ plus the area of the shaded part of the upper diagram in fig.~\ref{omarfig}, so $\bar {f}_0 (k_0)=k_0+a^2 < ak_0 + b$. 

Starting with the initial diagram $D_0$, we will modify the diagram according to the following algorithm. At each odd step we increase $k$, the length of the horizontal side of the rectangle, by one and add one to each part of the partition in the diagram (see fig.~\ref{omarfig}). At each even step we increase $k$ by one (in fig.~\ref{omarfig}). Let $D_i$ denote the diagram obtained after $i$ steps of the algorithm and $\bar {f}_{i}$ denote its complexity. After each odd step of the algorithm the shaded area in the diagram increases by $a-1$ (the complexity of the diagram increases by $a$), after each even step of the algorithm the shaded area in the diagram increases by $a$ (the complexity of the diagram increases by $a+1$). Thus $\bar {f}_{2i}(k_0+2i)=\bar {f}_0 (k_0) + i(2a+1)$. 

Since $\bar {f}_0(k_0) < ak_0 + b$, the integer $i_0=ak_0 + b- \bar {f}_0(k_0)$ is positive. The complexity of $D_{2i_0}$ is $\bar {f}_{2i_0}(k_0+2i_0)=\bar {f}_0 (k_0) + 2i_0a+i_0=\bar {f}_0 (k_0) + 2i_0a+ak_0 + b- \bar {f}_0(k_0)=a(k_0+2i_0)+b$. Let $K=k_0+2i_0$. The diagram $D_{2i_0}$ has complexity $\bar {f}_{2i_0}(K)=aK+b$. The diagram obtained by deleting the right-most column of $D_{2i_0}$ has the $(K-1)^{th}$ complexity equal to $\bar {f}_{2i_0} (K) -a$. 

The diagram $D_{2i_0}$ corresponds to a gap increasing prefix $v$ with the following property: any gap increasing word $w=vu$, such that $u$ has no subwords $10^n1$ with $n<K-1$, has subword complexity $f_w(K)=aK + b$ and $s_w(K-1)=a$. 

Next we will show that there exists $u$ such that $w=vu$ is gap increasing and
$f_w(n)=an + b$ for all $n \geq K$. By the construction above
$v=0^{n_0}10^{n_1}10^{n_2}1...10^{n_{2a-1}}1$ for some $n_0<n_1<...<n_{2a-1}$. 
Let $u=0^{n_{2a}}10^{n_{2a+1}}10^{n_{2a+2}}1...$, where 
$n_i=n_{i-a}+n_{i-a+1}+1$ for all $i \geq 2a$. Since $n_{2a+1}=n_a+n_{a+1}+1 \geq K-1$, $10^n1$ with $n<K-1$ is not a subword of $u$. Thus $w=vu$ has $K$th subword complexity $f_w(K)=aK + b$. By Lemma~\ref{lm:cond}, $s_w(n)=s_w(K-1)=a$ for all $n \geq K$. This implies that $f_w(n)=an + b$ for all $n \geq K$.

( $\leftarrow$ ) At last we have to show that there does not exist an infinite gap increasing word with subword complexity $f(n)=n+b$ ultimately.

Suppose there exists an infinite gap increasing word $w$ with subword complexity function $f_w$ and a positive integer $N$ such that $f_w(n)=n+b$ for $n>N$.   
For any integers $m$ and $n$ such that $n>m>N$, $f_w(n)-f_w(m)=n-m+\delta$, where $\delta$ is the difference between the number of distinct subwords of length $n$ and $m$ that contain at least two $1$'s. For $n$ large enough $\delta >0$. This proves the claim of the theorem.
\end{proof}

\section{Blockwise injective words}

\label{sec:block}

\begin{definition} A function $h: \N \rightarrow \N$ is called blockwise injective 
if $h(i)=h(j)$ for $i<j$ implies $h(i)=h(i+1)=...=h(j)$.
\end{definition}

Notice that the set of blockwise injective functions includes the set of injective 
functions and the set of non-decreasing functions. 

We will generalize the method of Proposition~\ref{pr:inj} to get a recurrence formula for
the subword complexity function of an infinite word $w$ whose gap function $g$ 
is blockwise injective. Here we will have to distinguish between unbounded and
bounded blockwise injective functions. To each unbounded blockwise injective function
$g:\N \rightarrow \N$ we assign two functions $j:\N \rightarrow \N$ and $p:\N
\rightarrow \N$ such that $j$ is the injective function which
assumes the same values and in the same order as $g$ does, and, 
for each $r \in \N$, $p(r)$ is the number of times $g$ assumes value $j(r)$. 
If $g:\N \rightarrow \N$ is a bounded blockwise injective function, let $b$ be the 
number of distinct values that $g$ assumes. Then  $j(r)$ is defined as above for 
$1 \leq r \leq b$, and $p(r)$ is defined for $1 \leq r \leq b-1$. 

\begin{proposition} Let $w$ be an infinite binary word whose gap function $g$
is blockwise injective. Let $j(r)$ and $p(r)$ be defined as above. Consider 
the following four systems of equations (with the convention that
$j(0)=0$):
\begin{equation}    
    j(r)\leq n < j(r)+ max \{ j(r-1),[sign (p(r)-1)]j(r) \},
    \label{eq:j1}
\end{equation}

\begin{equation}    
    \begin{cases} 
    n \geq j(r) \leq j(r')\\
    n < j(r)+max \{ j(r-1),[sign (p(r)-1)]j(r) \}  \\
    n < j(r)+max \{ j(r'-1),[sign (p(r')-1)]j(r') \} ,
    \end{cases}
    \label{eq:j1'}
\end{equation}

\begin{equation}    
    \begin{cases} 
    j(r+1)>j(r) \\
    2j(r)\leq n < p(r)j(r)+ min \{ j(r-1),j(r) \},
    \end{cases}
    \label{eq:j2}
\end{equation}
 
\begin{equation}    
    \begin{cases} 
    j(r+1)<j(r) \\
    j(r)+j(r+1)\leq n < [p(r)-1]j(r)+j(r+1)+ min \{ j(r-1),j(r) \}.
    \end{cases}
    \label{eq:j3}
\end{equation}

1) If $g$ is unbounded, then the number of right special factors of length $n$ of 
$w$ equals one plus the sum of the numbers of r solutions of equations~(\ref{eq:j1}), 
(\ref{eq:j2}), and (\ref{eq:j3}).

2) If $g$ is bounded, let $b$ be the number of distinct values that $g$ assumes and
$j_{max}=max \{j(r) \}$. Then the number of right special factors of length $n$ 
of $w$ is $I(n)+S_1(n)+S_2(n)$, where 
\begin{equation}
 I(n)=   
    \begin{cases} 
    1, \text{ if } n<j_{max}  \\
    0, \text{ if } n \geq j_{max},
    \end{cases}
\end{equation}
$S_1$ is the number of $r$'s for which there exists an $r' \geq 1$  such 
that eq.~(\ref{eq:j1'}) is satisfied, $S_2$ is the sum of the numbers of 
$r$ solutions of eq.~(\ref{eq:j2}) and (\ref{eq:j3}). 

\label{pr:bl}
\end{proposition}

\begin{proof}  1) First we consider the case when $g$ is unbounded. 

The word $0^n$ is a right special factor of $w$ for any integer 
$n \geq 1$ (this is the one which is not counted by the solutions of the system). 

We will show that the number of right special factors of length $n$ of $w$ which 
contain exactly one $1$ equals the number of integer $r$ solutions of eq.~(\ref{eq:j1}). 
A binary word which contains exactly one $1$ is a right special factor 
of $w$ if and only if it has the form $0^k10^{j(r)-1}$, where $r \geq 1$ and $k$ 
satisfies one of the following: 

(a) $p(r)=1$ and $0 \leq k < j(r-1)$.

(b) $p(r) > 1$ and $0 \leq k < max \{ j(r-1),j(r)\}$. (Reminder: $j(0)=0$) 

The number of subwords of $w$ of length $n$ of form $0^k10^{j(r)-1}$ 
which satisfy conditions (a) and (b) equals the number of $r$'s 
which satisfy equations~(\ref{eq:one11}) and (\ref{eq:one12}) respectively.

\begin{equation}    
    \begin{cases} 
    p(r)=1 \\
    0 \leq n-j(r) < j(r-1),
    \end{cases}
    \label{eq:one11}
\end{equation}
\begin{equation}    
    \begin{cases} 
    p(r)>1  \\
    0 \leq n-j(r) < max \{ j(r-1),j(r) \} .
    \end{cases}
    \label{eq:one12}
\end{equation}
      
Equations~(\ref{eq:one11}) and (\ref{eq:one12}) have disjoint sets of $r$ solutions 
and eq.~(\ref{eq:j1}) combines the solutions of both, which explains why the number of
solutions of eq.~(\ref{eq:j1}) is the same as the number of right special factors of 
length $n$ of $w$ which contain exactly one $1$.

Next we shall count the number of right special factors of length $n$ of $w$ which 
contain at least two $1$'s. Here we use the convention that if $u$ is a finite word, 
then $u^0= \epsilon$ (the empty word).

If a right special factor of $w$ contains at least three $1$'s, the number of $0$'s
between every two consecutive $1$'s in this right special factor should be the same. 
Indeed, since $g$ is blockwise injective, any subword $v$ of $w$ which has 
$10^s10^t1$ ($s \ne t$) as a subword, can occur only once in $w$, thus $v$ 
cannot be a right special factor. Therefore a right special factor of $w$ which 
contains at least two $1$'s is necessarily, but not sufficiently, 
of form~(\ref{type1}) or (\ref{type2}).
\begin{equation} 
0^k(10^{j(r)-1})^m,\hspace{0.2in} k \geq 0,\hspace{0.2in} r \geq 1,\hspace{0.2in} 
2 \leq m \leq p(r);
\label{type1}
\end{equation}
\begin{equation}
0^k(10^{j(r)-1})^m10^{j(r+1)-1},\hspace{0.2in}  k \geq 0,\hspace{0.2in} r \geq 1,
\hspace{0.2in} 1 \leq m \leq p(r).
\label{type2}
\end{equation}

We will show that the number of right special factors of length $n$ of
form~(\ref{type1})  equals the number of $r$ solutions of eq.~(\ref{eq:j2}),
and the number of right  special factors of length $n$ of form~(\ref{type2})
equals the number of $r$  solutions of eq.~(\ref{eq:j3}).

First we shall count the number of right special factors of length $n$ and 
form~(\ref{type1}). Let $v=0^k(10^{j(r)-1})^m$ be a subword of $w$, where $k
\geq 0$, $r \geq 1$, $2 \leq m \leq p(r)$. The subword $v$ is followed by $1$
in its leftmost occurrence in $w$. It is  followed by $0$ (which can happen
only in its rightmost occurrence) if and only if
\begin{equation}
j(r+1)>j(r). 
\label{eq:b1}
\end{equation}

Hence a binary word of form~(\ref{type1}) is a special factor of $w$ if and
only if it is a subword of $w$ and (\ref{eq:b1}) holds. For a fixed $r$ there
exists a subword of $w$ of form~(\ref{type1}) and length $n$, if and only if 
\begin{equation}
2j(r) \leq n < p(r)j(r)+ min \{ j(r-1),j(r) \}. 
\label{eq:b2}
\end{equation}

It should be mentioned, that for every $r$ that satisfies eq.~(\ref{eq:b2}),
there exists a unique subword of $w$ of form~(\ref{type1}) and length $n$. 

Notice that system~(\ref{eq:j2}) is a combination of equations~(\ref{eq:b1})
and (\ref{eq:b2}), and there is a bijective correspondence
between right special factors of length $n$ and form~(\ref{type1}) and $r$
solutions of system~(\ref{eq:j2}). 

At last we shall count the number of right special factors of length $n$ of 
form~(\ref{type2}). Let $u=0^k(10^{j(r)-1})^m10^{j(r+1)-1}$ be a subword of
$w$, where $k \geq 0$, $r \geq 1$, $1 \leq m \leq p(r)$. The subword $u$ is
followed by $1$ in its rightmost occurrence in $w$. It is  followed by $0$
(which can happen only in its leftmost occurrence) if and only if
\begin{equation}
j(r+1)<j(r). 
\label{eq:c1}
\end{equation} 

Hence a subword of $w$ of form~(\ref{type2}) is a right special factor if and 
only if (\ref{eq:c1}) holds. For a fixed $r$ there exists a subword of $w$ of
form~(\ref{type2}) and length $n$ (and there can be only one such subword) 
if and only if 

\begin{equation}
j(r)+j(r+1)\leq n < [p(r)-1]j(r)+j(r+1)+ min \{ j(r-1),j(r) \}.
\label{eq:c2}
\end{equation}

Thus for every $r$ there exists a right special factor of form~(\ref{type2}) 
and length $n$ (this special factor happens to be unique) if and only if 
$r$ satisfies system~(\ref{eq:j3}). That proves that the right special 
factors of $w$ of length $n$ are counted by the number of $r$ solutions of
system~(\ref{eq:j3}).

\vspace{.1in}

2) Next we consider the case when $g$ is bounded. 

It is clear that $0^n$ is a right special factor of $w$ if and 
only if $n<j_{max}$, this accounts for $I(n)$. The number of right special 
factors of length $n$ which contain at least two $1$'s equals $S_2$, the 
argument is the same as in the case when $g$ is unbounded.

At last we have to show that the number of right special factors of length 
$n$ which contain exactly one $1$ equals the number of $r$'s for which there exists
an $r' \geq 1$ such that eq.~(\ref{eq:j1'}) is satisfied.

Let $v$ be a binary word of length $n$ which contains exactly one $1$, then $v$
is a right special factor of $w$ 
if and only if $v=0^{n-j(r)}10^{j(r)-1}$ and $v0=0^{n-j(r)}10^{j(r)}$ is a subword
of $w$, that is
\begin{equation}
0\leq n-j(r) < max \{ j(r-1),[sign (p(r)-1)]j(r) \}. 
\label{eq:bl21}
\end{equation}
and there is $r' \geq 1$ such that 
\begin{equation}    
    \begin{cases} 
    j(r') \geq j(r)\\
    n-j(r) < max \{ j(r'-1),[sign (p(r')-1)]j(r') \} .
    \end{cases}
    \label{eq:bl22}
\end{equation}

Thus the number of right special factors of $w$ which contain exactly one $1$
equals the number of $r$ solutions of eq.~(\ref{eq:bl21}) for which there is $r' \geq
1$ such that eq.~(\ref{eq:bl22}) is satisfied and, because eq.~(\ref{eq:j1'})
is obtained by combining eq.~(\ref{eq:bl21}) and eq.~(\ref{eq:bl21}), the
claim is proved.  
\end{proof}

\begin{corollary} Let $w$ be an infinite binary word whose gap function $g$
is non-decreasing and $j(r)$, $p(r)$ (and $b$ if $g$ is bounded) be 
defined as before. Let $s_w(n)$ denote the number of right special factors of 
length $n$ of $w$.

If $g$ is unbounded, then, for every natural $n$, $s_w(n)$ equals one plus the
number of integer $r$ solutions of   
\begin{equation}    
    \begin{cases} 
    j(r)\leq n \\
    j(r-1)+p(r)j(r)\geq n+1 .
    \end{cases}
    \label{eq:nondec1}
\end{equation}
If $g$ is bounded, for every natural number $n$ consider the inequality system 
\begin{equation}
    \begin{cases} 
    r<b \\
    j(r)\leq n \\
    j(r-1)+p(r)j(r)\geq n+1 .
    \end{cases}
    \label{eq:nondec2}
\end{equation}
If $1 \leq n \leq j(b)-1$, then $s_w(n)$ equals one plus the number of 
integer $r$ solutions of eq.~(\ref{eq:nondec2}). If $n \geq j(b)$, then $s_w(n)$ is 
just the number of integer solutions of eq.~(\ref{eq:nondec2}).
\label{cor:nondec}
\end{corollary}
\begin{proof} First consider the case when $g$ is unbounded. By Proposition~\ref{pr:bl}
$s_w(n)$ equals one plus the number of $r$'s which satisfy one of the
following:
\begin{equation}    
    \begin{cases} 
    p(r)=1 \\
    0 \leq n< j(r-1)+j(r),
    \end{cases}
    \label{eq:nondec3} 
\end{equation}

\begin{equation}    
    \begin{cases} 
    p(r)>1  \\
    0 \leq n< 2j(r),
    \end{cases}
    \label{eq:nondec4}
\end{equation}

\begin{equation}
2j(r)\leq n < p(r)j(r)+j(r-1).
\label{eq:nondec5}
\end{equation}

The set of $r$ solutions of system~(\ref{eq:nondec1}) is the disjoint union of
the sets of $r$ solutions of equations~(\ref{eq:nondec3}), (\ref{eq:nondec4}) 
and (\ref{eq:nondec5}).

Next we consider the case when $g$ is bounded. Because $g$ is
non-decreasing, if $r$ satisfies eq.~(\ref{eq:j1}), then $r$ and  
$r'=r+1$ satisfy eq.~(\ref{eq:j1'}). Thus the set of $r$'s which satisfy of
eq.~(\ref{eq:j1'}) is the same set that satisfies eq.~(\ref{eq:j1}).
By  Proposition~\ref{pr:bl}, the number of
right special factors of $w$ of length $n$ equals the number 
of $r$ solutions of eq.~(\ref{eq:nondec1}) with the restriction $r<b$ 
(because p(r) is defined only for $r<b$) plus $I(n)$, where  
\begin{equation*}
  I(n)=   
    \begin{cases}      
    1, \text{ if } n \leq j(b)-1  \\     
    0, \text{ if } n > j(b).     
\end{cases} 
\end{equation*} 
\end{proof}

\begin{example} Consider the infinite binary word whose non-decreasing gap 
function is given by $j(r)=r+1$ and $p(r)=r$. 
$$w=0(10^2)^2(10^3)^3(10^4)^4...$$

By Corollary~\ref{cor:nondec} the number of right special factors of
length $n$ of $w$ equals the number of solutions of the system
\begin{equation*}    
    \begin{cases} 
    r+1\leq n \\
    r+r(r+1)\geq n+1 .
    \end{cases}
\end{equation*}

This system has $n- \left \lceil \sqrt {n-2} \right \rceil$ solutions, thus 
$$f_w(n+1)=f_w(n)+n- \left \lceil \sqrt {n-2} \right \rceil.$$
\end{example}

\begin{proposition} Let $w$ be an infinite binary word with gap function 
$g$  and subword complexity function $f$. If $g$ is non-decreasing, then 
$0 \leq \Delta f(n) \leq n+1.$ If $g$ is blockwise injective, then 
$0 \leq \Delta f(n) \leq n(n-1)/2 + 1.$ 
\label{pr:simnec2}   
\end{proposition}

\begin{proof} The proof is similar to the proof of Proposition~\ref{pr:simnec}.
$\Delta f(n)$ is the number of right special factors of length $n$ of $w$.

If $g$ is non-decreasing, $\Delta f(n)$ is at least zero ($\Delta f(n)=0$ for some $n$ if and only if $w$ is ultimately periodic, which happens if and only if $g$ is bounded). Also $\Delta f(n)$ is at most the number of subwords of $w$ of 
form $0^k(10^l)^m1$, where $k,l,m \geq 0$, $k \leq l$ and $k+(l+1)m=n$ (this follows from the proof of Proposition~\ref{pr:bl}). There are at most $n+1$ such subwords because $l$ completely determines the subword $0^k(10^l)^m1$ and there are $n+1$ choices for $l$, hence $0 \leq \Delta f(n) \leq n+1.$

If $g$ is blockwise injective, $\Delta f(n)$ is at least zero and at most the number of subwords of $w$ of form $0^k(10^l)^m(10^r)^i1$, 
where $k,l,m,r \geq 0$, $l \ne r$, $k+(l+1)m+(r+1)i=n$, $i=0$ if $k=n$, and $i=1$ otherwise (again follows from the proof of Proposition~\ref{pr:bl}). 
There are at most $n(n+1)/2+1$ such subwords because $l$ and $r$ determine 
the subword $0^k(10^l)^m(10^r)^i1$ completely (except for the case when $k=n$) and 
$0 \leq l+r \leq n-1$, thus there are at most $n(n+1)/2$ such pairs of $l$ and 
$r$ (add one for the case $k=n$). Hence $0 \leq \Delta f(n) \leq n(n+1)/2+1.$
\end{proof}

\begin{corollary} Let $w$ be an infinite binary word with a blockwise injective gap function and subword complexity function $f$. Then $1 \leq f(n) \leq \frac {n^3}{6} + \frac {5n}{6} +1$.
\end{corollary}

\begin{proposition} Let $a$ and $b$ be two integers. There exists an infinite binary word $w$ with a blockwise injective gap function and subword complexity function $f_w(n)=an+b$ ultimately if and only if one of the two following conditions holds:

1) $a \geq 2$. 

2) $a=0$ and $b \geq 1$.
\end{proposition}
\begin{proof}
It was proved in Theorem~\ref{th:anb} that for any integers $a \geq 2$ and $b$ there exists a gap increasing word $w$ with subword complexity $an+b$ ultimately. 
Since $w$ is gap increasing, the subword complexity function of $w$ is blockwise injective.

For any integer $b \geq 2$ the binary infinite word $w_b=(10^{b-1})^{\infty}$ has the subword complexity $f(n)=b$ for all $n \geq b-1$. Clearly the gap function of $w_b$ is blockwise injective.   

At last we have to show that an infinite binary word, whose gap function is blockwise injective, cannot have subword complexity $f(n)=n+b$ ultimately. An 
infinite word $w$ whose gap function is blockwise injective is either ultimately 
periodic or $0^i10^j$ is a subword of $w$ for all $i$ and $j$. If $w$ is ultimately periodic, then its subword complexity function is ultimately constant and the claim is proved. If $0^i10^j$ is a subword of $w$ for all $i$ and $j$, then $w$ contains exactly $n+1$ distinct subwords of length $n$ that contain at most one $1$. Thus for any $m<n$, $f_w(n)-f_w(m)=n-m+\delta_{m,n}$, where $\delta_{m,n}$ is the difference between the number of distinct subwords of length $n$ and $m$ that contain at least two $1$'s. To every subword $v$ of $w$ of length $m$ that contains at least two $1$'s we can put in correspondence the subword $u$ of $w$ of length $n$ with prefix $v$ (this $u$ contains at least two $1$'s). For any $m$, there exists $n$ large enough such that $10^{n-2}1$ is a subword of $w$ whose prefix of length $m$ contains only one $1$, hence $\delta_{m,n} > 1$. This means that there does not exist $m$ such that for all $n \geq m$ the subword complexity of $w$ is $f(n)=n+b$. 
\end{proof}

\end{document}